\documentclass[a4paper,12pt]{article} 
\usepackage{amsmath,amsfonts,bbm} 
\usepackage {ucs}
\usepackage[utf8x]{inputenc}

\newtheorem{thm}{Theorem}[section]   
\newtheorem{cor}[thm]{Corollary}   
\newtheorem{prop}[thm]{Proposition}   
\newtheorem{lm}[thm]{Lemma}

\newcommand{\RR}{\mathbb{R}}   
\newcommand{\di}{\displaystyle}

\newcommand{\e}{\varepsilon}
\def\un{{\mathbf{1}}}

\begin{document}   
   
\title{\textbf{Sharp large time behaviour in $N$-dimensional Fisher-KPP equations}}
\author{
{\bf Jean-Michel Roquejoffre}\\
Institut de Math\'ematiques de Toulouse; UMR 5219\\
Universit\'e de Toulouse; CNRS\\ 
Universit\'e Toulouse III,
118 route de Narbonne, 31062 Toulouse, France \\ 
\texttt{jean-michel.roquejoffre@math.univ-toulouse.fr}\\ 
\\[.5mm]
{\bf Luca Rossi} \\ 
Centre d'Analyse et de Math\'ematique Sociales; UMR 8557\\
Paris Sciences et Lettres; CNRS\\ 
EHESS,
54 Bv. Raspail, 75006 Paris, France\\ 
\texttt{luca.rossi@ehess.fr}\\
\\[.5mm]
{\bf Violaine Roussier-Michon} \\ 
Institut de Math\'ematiques de Toulouse; UMR 5219\\
Universit\'e de Toulouse; CNRS\\ 
INSA Toulouse,
135 av. Rangueil, 31077 Toulouse, France\\ 
\texttt{roussier@insa-toulouse.fr}
}
\date{{\it Dedicated to L. Caffarelli, as a sign of friendship, admiration and respect.}
}
\maketitle 
\begin{abstract}
We study the large time behaviour of the Fisher-KPP equation $\partial_t u=\Delta u +u-u^2$ in spatial dimension $N$, when the initial datum is compactly supported. We prove the existence of a Lipschitz function $s$ of the unit sphere, such that $u(t,x)$ converges, as $t$ goes to infinity, to
 $U_{c_*}\bigg(|x|-c_*t +\di\frac{N+2}{c_*} \mathrm{ln}t + s^\infty\Big(\frac{x}{|x|}\Big)\bigg)$, where $U_{c*}$ is the 1D travelling front with minimal speed $c_*=2$. This extends an earlier result of G\"artner.
\end{abstract}

\section{Introduction}    
\label{s1}

\noindent
The paper is devoted to the large time behaviour of the solution of the reaction-diffusion equation 
\begin{align}
\label{KPP} \partial_t u= \Delta u +f(u), & \quad t>1 \, , \,  x \in \RR^N\\
\notag u(1,x)=u_0(x), & \quad \quad \quad  \quad x \in \RR^N
\end{align}
We will take
$$
f(u)=u(1-u);
$$
thus  $f$ is, in reference to the pioneering paper \cite{KPP}, said to be of the Fisher-KPP type. The initial datum $u_0$ is in ${\cal C}(\RR^N)$ and there exist $0<R_1<R_2$ such that
\begin{equation}
\label{assumption}
\forall x \in \RR^N \, , \quad \un_{B_{R_1}}(x)\leq u_0(x) \leq \un_{B_{R_2}}(x),
\end{equation}
where $\un_A$ is the indicator of the set $A$ and $B_{R}$ the ball of $\RR^N$ of radius $R$. 
By the maximum principle and the standard theory of parabolic equations (see for instance \cite{Hy}), equation \eqref{KPP} has a unique classical solution $u(t,x)$ in $ {\cal C}([1,+\infty[ \times \RR^N, (0,1))$ emanating from $u_0$. The first and most general result is due to Aronson and Weinberger \cite{aw}. The solution $u$ spreads at the speed $c^*=2\sqrt{f'(0)}=2$ in the sense that 
$$
\min_{|x|\leq ct} u(t,x) \to 1 \mbox{ as } t \to +\infty \, , \mbox{ for all } 0 \leq c < c^*
$$
and 
$$
\sup_{|x|\geq ct} u(t,x) \to 0 \mbox{ as } t \to +\infty \, , \mbox{ for all }  c > c^*.
$$
The goal of this paper is to sharpen this result.

\noindent Let us briefly recall what happens as time goes to infinity when $N=1$.  Equation \eqref{KPP} with $N=1$ reads
\begin{equation}
\label{KPP1D}
\partial_t u=\partial_{xx} u +f(u),\quad t>1\, , \, x\in\RR.
\end{equation}
It admits one-dimensional travelling fronts $U(x-ct)$ if and only if $c\geq c^*=2$ where the profile~$U$, depending on $c$, satisfies
\begin{equation}
\label{onde}
U'' +c \,U'+f(U)=0 , \quad x \in  \RR,
\end{equation}
together with the conditions at infinity 
\begin{equation} 
\label{CL onde}
\lim\limits_{x \to - \infty} U(x)=1 \mbox{ and } 
\lim\limits_{x \to + \infty} U(x)=0. 
\end{equation} 
Any solution $U$ to \eqref{onde}-\eqref{CL onde} is a shift of a fixed profile $U_c$: $U(x)=U_c(x+s)$ with some fixed $s \in \RR$. The profile $U_{c^*}$ at minimal speed $c^*=2$ satisfies, up to translation,
$$
U_{c^*}(x)=(x+K) \,  e^{- x} +O(e^{-(1+\delta_0)x})\, , \, \mbox{ as } x \to +\infty
$$
for some universal constants $K \in \RR$ and  $\delta_0>0$.
The large time behaviour of  \eqref{KPP1D} has a history of important contributions, we are only going to list two lasting ones.
The  first  is  the  paper of Kolmogorov, Petrovskii and Piskunov \cite{KPP}. They proved 
that the solution of~\eqref{KPP1D} starting from the initial datum
$\mathbbm{1}_{(-\infty,0]}$ converges to $U_{c*}$ in shape: there is a   function
$$
\sigma^\infty(t)=2t+o_{t\to+\infty}(t),
$$
such that 
\[
\di
\lim_{t\to+\infty}u(t,x + \sigma^\infty(t))=U_{c^*}(x) \quad\mbox{  uniformly in }x \in \RR.
\]
The second contribution makes precise the $\sigma^\infty(t)$: in \cite{Bramson2},   Bramson
 proves the existence of a constant $x_\infty$, depending on $u_0$, such that
\begin{equation}
\label{e1.000}
\sigma^\infty(t)=2t-\di
\frac32\ln t-x_\infty+o_{t\to+\infty}(1).
\end{equation}
Formula \eqref{e1.000} was proved through elaborate probabilistic arguments.  As said before, the problem, as well as more complex variants of it, are currently the subject of intense investigations. See for instance \cite{NRR-Brezis} for an account of them.

In several space dimensions, the asymptotics have been pushed less far. 
In the framework of the Fisher-KPP equation that we are studying, the
Aronson-Weinberger result is made precise  up to $O(1)$ terms in G\"artner \cite{Ga}. If $N$ is the space dimension, the main result of \cite{Ga} is that, for every $\lambda\in(0,1)$, the level set $\{u=\lambda\}$ is trapped, for large times, between two spheres of radius 
$$
R(t)=c^*t-\frac{N+2}{c^*}{\mathrm{ln}}t+O_{t\to+\infty}(1).
$$ 
The $O_{t\to+\infty}(1)$ terms are not studied. G\"artner's contribution is probabilistic, and a PDE proof of his result is provided by Ducrot \cite{D}, adapting to higher dimension 
the proof of (a weaker version of) Bramson's formula \eqref{e1.000}, given by F. Hamel, J. Nolen, L. Ryzhik and the first author in \cite{HNRR1}. 

\noindent When the coefficients of the equation actually depend on $x$ 
in a periodic fashion, as for instance for the equation
$$
\partial_t u=\Delta u+\mu(x)u-u^2,\quad t>0,\ x\in\RR^N,
$$
with $\mu$ periodic and positive (actually, more general assumptions on $\mu$ can be allowed, as well as inhomogeneous diffusion terms, or the presence of advection), a lot is now known on the spreading speed, or, in other words, the position of the level sets up to $O_{t\to+\infty}(1)$ terms. The first result in this direction is Freidlin-G\"artner \cite{FG}, which gives, through a probabilistic approach, an almost explicit expression (the Freidlin-G\"artner formula) of the spreading speed in each direction. Several proofs and generalisations of this formula have been given, by various approaches: viscosity solutions \cite{EvS}, abstract dynamical systems \cite{W}, PDE approach \cite{BH}, \cite{Ross2}. Let us mention an important contribution~\cite{Sha}, which generalises G\"artner's result to periodic functions $\mu(x)$, by computing the relevant logarithmic shift. This work also generalises \cite{HNRR2}, a contribution that computes the shift for periodic $\mu$, but in one space dimension.

\noindent Coming back to \eqref{KPP}, the goal of the present paper is to make precise the $O_{t\to+\infty}(1)$ in G\"artner's expansion. Our result is the
\begin{thm}
\label{thm1}
Let $u_0$ satisfy assumption \eqref{assumption}. There is a Lipschitz function $s^\infty(\Theta)$, defined on the unit sphere of $\RR^N$, such that the solution $u$ of \eqref{KPP} emanating from $u_0$ satisfies
$$
u(t,x)=U_{c_*}\biggl(\vert x\vert-c_*t+\frac{N+2}{c_*}{\mathrm{ln}}t+
s^\infty\Big(\frac{x}{\vert x\vert}\Big)\biggl)+o_{t\to+\infty}(1),
$$
with $c_*=2$.
\end{thm}

\noindent  This completes the result of \cite{Ga}. At this stage, let us anticipate the proof of the theorem, and let us give a brief explanation of the logarithmic shift observed here: it can be decomposed into two shifts having different origins. One is due to the curvature term $\frac{N-1}{c_*} \ln t$, it systematically arises in this type of large time issues for reaction-diffusion equations, the nonlinearity $f$ does not need to be of the KPP type. See for instance \cite{VRM}, \cite{yagisita}. The other is the one-dimensional shift $\frac3{c_*}{\mathrm{ln}}t$, it is typical of the KPP nonlinearity. All this will be made clearer in Section 2.

\noindent Theorem \ref{thm1} is in contrast with a recent paper \cite {JMR_VRM2} of the first and third authors, which studies \eqref{KPP} when the initial datum is trapped between two planar travelling waves. In this setting, the logarithmic shift is $\di\frac32{\mathrm{ln}}t$, as in the one-dimensional case. However, the dynamics beyond the logarithmic shift is given by that of the heat equation on the whole line. This last equation, though extremely well-behaved as far as the regularity of its solutions is concerned, exhibits solutions that do not converge, as time goes to infinity, to anything. However, this last feature holds for reaction-diffusion that need not be of the KPP type, see \cite{JMR-VRM}.

\noindent Before starting the proof of our results, let us mention that it would be certainly interesting to understand sharper asymptotics of $u(t,x)$. In one space dimension, a full expansion has been proposed, in the formal style, in \cite{ES}, or with another approach in \cite{BBD}. The next term in the expansion of the shift is computed, in a mathematically rigorous way, in \cite{NRR2}. The expansion is pushed even further in \cite{Gra}.

\noindent Let us also say that the observed behaviour is quite typical of Fisher-KPP equations with second-order linear diffusion. Another important class of nonlinearities $f(u)$ in \eqref{KPP} satisfies $f(0)=f(1)=0$, $f'(0)<0$, $f'(1)<0$, with $\di\int_0^1f(u)du>0$. A typical example  is
$$
f(u)=u(u-\theta)(1-u),\quad 0<\theta<\frac12.
$$
A statement of the same type as Theorem \ref{thm1} is  \cite{VRM}, with the important difference that the logarithmic delay is solely due to the curvature terms; the dynamics beyond the shift is the same as the one presented in Theorem \ref{thm1}. And, although the phenomenon does not look so remote to the one displayed in \cite{VRM}, it is quite different in nature, as the convergence to the wave is dictated by what happens in the region where the solution takes intermediate values. A similar, and recent contribution \cite{DQZ} treats the porous medium equation with Fisher-KPP nonlinearity; although the nonlinearity is the same as in the present paper, the result is of the type of \cite{VRM} (although the dynamics beyond the shift is not made precise when the initial datum is nonradial), this is due to the fact that the solution does not have a tail that would govern the overall dynamics. We end this series of remarks by recalling a result of Jones \cite{J}, stating that 
the level sets of the solution of \eqref{KPP}, whatever the nonlinearity is, will have oscillations only of the size $O_{t\to+\infty}(1)$. This is a consequence of the following fact: if $\lambda$ is a regular value of $u$, the normal to the $\lambda$-level set of $u$ meets the support of the initial datum. A very simple proof of this fact is given by Berestycki in \cite{B}.

\noindent In the next section, we transform the equations so as to uncover the basic mechanism at work, namely the fact that the whole phenomenon is dictated by the tail of the solution. The subsequent sections are different steps of the proof of Theorem \ref{thm1}, this will be explained in more detail in Section 2.

\medskip
\noindent
{\bf Acknowledgement.} JMR and LR are supported by  the European Union's Seventh 
Framework Programme (FP/2007-2013) / ERC Grant
Agreement n. 321186 - ReaDi - ``Reaction-Diffusion Equations, Propagation and Modelling''. VRM is supported by the ANR project NONLOCAL ANR-14-CE25-0013.

\section{Preparation of the equations, method of proof}
\label{s2}

There is a sequence of transformations that bring the equations under the \eqref{KPP} to 
a form that will make evident that the region $|x|\sim\sqrt t$ in the moving frame, that we will subsequently call the diffusive zone, dictates the hole dynamics. 
\begin{enumerate}
\item We first use the polar coordinates 
$$
x\mapsto (r=\vert x\vert>0,\Theta=\frac{x}{\vert x\vert}\in \mathbb{S}^{N-1})
$$
then \eqref{KPP} becomes
$$
\partial_t u=\partial_{rr}u +\frac{N-1}r\partial_ru +\frac{\Delta_\Theta u}{r^2}+ u-u^2,\quad{t>1,r>0,\Theta\in \mathbb{S}^{N-1}}.
$$
Here, $\Delta_\Theta$ is the Laplace-Beltrami operator on the unit sphere of $\RR^N$. Its precise expression will not be needed in the sequel.
\item Let us believe that the transition zone where $u$ decreases from 1 to 0 is located around $R(t)=2t- k \ln t $ ($k$ to be chosen later) and choose the change of variables $r'=r-R(t)$ and $u(t,r,\Theta)=u_1(t,r-R(t),\Theta)$. We drop the primes and indexes, and \eqref{KPP} becomes
\begin{equation}
\label{ref mvt} 
\partial_t u= \partial_{rr}u+\frac{N-1}{r+2t-k{\mathrm{ln}}t}\partial_ru+(2-\frac{k}t)\partial_r u+\frac{\Delta_\Theta u}{(r+2t-k{\mathrm{ln}}t)^2}+u-u^2.
\end{equation}
The equation is valid for $t>1$, $r>-2t+k{\mathrm{ln}}t$, and $\Theta\in \mathbb{S}^{N-1}$.
 \item As is by now classical, we take  out the exponential decay of the wave $U_{c^*}$, and set $u(t,r,\Theta)=e^{-r} v(t,r,\Theta)$; \eqref{ref mvt} thus becomes
\begin{equation}
\label{dec onde} 
\partial_t v= \partial_{rr} v + (\frac{N-1}{r+2t-k{\mathrm{ln}}t}-\frac{k}t) \left(\partial_r v - v\right)+\frac{\Delta_\Theta v}{(r+2t-k{\mathrm{ln}}t)^2} -e^{-r} v^2,
\end{equation}
with initial datum $v(1,r,\Theta)=e^r u_0(r+2,\Theta)$.
\item We now choose $k$. Our first guess is that the term in $\Delta_\Theta v$ will not matter too much, because it decays like $t^{-2}$ (an integrable power of $t$), except in the zone $r\sim -2t$, where we know (for instance \cite{aw}) that $u(t,r,\Theta)$ goes to 1 as $t\to+\infty$. Hence we expect the dynamics to be like that of the one-dimensional equation. On the other hand, in the advection term, the quantity
$
\di\frac{N-1}{r+2t-k{\mathrm{ln}}t}
$
is nonintegrable in $t$, except for extremely large $r$. Thus we wish to balance it with the $\di\frac{k}t$ term. However, instructed by the large time behaviour in one space dimension, we keep in mind that we should keep the quantity $-\di\frac3{2t}$ factoring $\partial_rv-v$. Hence we choose
\begin{equation}
\label{e2.10000}
\frac{N-1}2-k=-\frac32,
\end{equation}
hence 
\begin{equation}
\label{e2.1}
k=\di\frac{N+2}2=\frac{N+2}{c_*}.
\end{equation}
In the sequel, we will keep the notation $k$, keeping in mind that $k$ is defined by \eqref{e2.1}.
\item Finally, in order to study \eqref{dec onde} in the diffusive zone, that is, the region $r\sim\sqrt t$,
 we use the self-similar variables 
$\displaystyle \xi=\frac{r}{\sqrt{t}}$, $\tau=\ln t $. The variable $\Theta$ is unchanged:
\begin{equation}
\label{var auto sim}
w(\tau,\xi,\Theta)=w \left(\ln t, \frac{r}{\sqrt{t}}, \Theta \right)=\frac{1}{\sqrt{t}} v(t,r,\Theta)
\end{equation}
Then \eqref{dec onde} becomes 
\begin{equation}
\label{shift} 
 \partial_{\tau} w+{\cal L} w = \frac{e^{\tau}\Delta_\Theta w}{(2e^\tau+\xi e^{\tau/2}-k\tau)^2}+h(\tau,\xi) e^{-\frac{\tau}{2}} \partial_\xi w-\biggl(h(\tau,\xi)+\frac32\biggl)w - e^{\frac{3}{2}\tau - \xi e^{\frac{\tau}{2}}} w^2,
\end{equation}
where 
$${\cal L} w = -\partial_{\xi\xi} w -\frac{\xi}{2} \partial_{\xi} w - w.$$ 
The function $h$ can easily be computed, although its expression is lengthy. It satisfies, for all $\delta\in(0,1/2)$:
\begin{equation}
\label{e2.1000}
h(\tau,\xi)=
\left\{
\begin{array}{rll}
&-\di\frac32+O(e^{-\delta\tau})\ \hbox{for $\xi\leq e^{(\frac12-\delta)\tau}$, that is, $r\leq t^{1-\delta}$,}\\
&O(1)\ \hbox{for $\xi\geq e^{(\frac12-\delta)\tau}$, that is, $r\geq t^{1-\delta}$}
\end{array}
\right.
\end{equation}
Of course the range of $\xi$ should be restricted in the negative direction, that is, $\xi > -2e^{\frac{\tau}{2}}+k\tau e^{-\frac{\tau}{2}}$, thus a very negative quantity if $\tau$ is very large. As the range of negative $\xi$ that will preoccupy us will be extremely modest (we will always have $\xi\geq-e^{-(\frac12-\delta)\tau}$, that is, $r\geq -t^{\delta}$) we will not mention this constraint in the sequel. Finally, the initial datum at $\tau=0$  is
$$
w_0(\xi,\Theta)=e^{\xi} u_0(\xi +2,\Theta).
$$
\end{enumerate}

Let us say a word about the strategy of the proof of Theorem \ref{thm1}. It will be inspired from the ideas of \cite{NRR-Brezis} in one space dimension, with some actual novelties due to the transverse variable. Our main step will be to prove the existence of a Lipschitz function $\alpha^\infty(\Theta)$ such that 
$$
w(\tau,\xi,\Theta)\longrightarrow_{\tau\to+\infty} \alpha^\infty(\Theta)\xi^+e^{-\xi^2/4},\quad\hbox{in $\{\xi\geq e^{-(\frac12-\delta)\tau}\}$},
$$
where $\delta>0$ is arbitrarily small. We will see that we cannot expect a better regularity than Lipschitz. The parallel step in \cite{NRR-Brezis} for $N=1$ was to prove, for the equation
 $$\partial_{\tau} w + {\cal L} w= -\frac{3}{2} e^{-\frac{\tau}{2}} \partial_{\xi} w - e^{\frac{3}{2}\tau - \xi e^{\frac{\tau}{2}}} w^2 \, , \quad  \tau>0 \, ,\quad \xi\in \RR,$$
 the existence of a constant $\alpha^\infty>0$ such that
$$
w(\tau,\xi)\longrightarrow_{\tau\to+\infty} \alpha^\infty\xi^+e^{-\xi^2/4},\quad\hbox{in $\{\xi\geq e^{-(\frac12-\delta)\tau}\}$}.
$$
The main effort was to prove the compactness of the trajectories $(w(\tau+T,\xi))_{T>0}$ as $T\to+\infty$; because the limiting trajectories satisfied the Dirichlet heat equation in self-similar variables, this entailed the convergence to a single Gaussian. To prove the compactness, we used a pair of sub/super solutions very much in the spirit of Fife-McLeod \cite{FML}; that one could actually use ideas from the analysis of bistable equations came as a surprise to us. In \cite{NRR-Brezis}, the barriers that we devised, a sort of miracle in the sign of the disturbances (that is, the exponential correction in the function $h$) allowed them to be sub and super solutions all the way down to $\xi=0$. Because we are dealing with a more complex equation, we do not want to rely on sign considerations, and we devise a pair of barriers that are sub and super solutions for more robust reasons than in \cite{NRR-Brezis}. They rely on a technical innovation in the vicinity of $\xi=0$, that is, if one thinks very much about the Fife-McLeod sub/super solutions, quite in the spirit of \cite{FML} once again. Once this is constructed, an additional issue will be to deal with the variable $\Theta$: as $\tau\to+\infty$, the Laplace-Beltrami operator will disappear from the asymptotic equations. That is, asymptotic regularity in $\Theta$ will have to be retrieved with bare hands, and this is why we cannot expect much more than Lipschitz regularity in $\Theta$.   
 
\noindent Once convergence in the diffusive area is under control, the next step is to fix the translation $\sigma^\infty(t,\Theta)$. We choose it  such that 
$$
U_{c_*}(r+\sigma^\infty(t,\Theta))\biggl\vert_{r=t^\delta}=e^{-r}v(t,r,\Theta)\biggl\vert_{r=t^\delta}.
$$
That is, 
$$
\sigma^\infty(t,\Theta)=-{\mathrm{ln}}\alpha^\infty(\Theta)+O(t^{-\delta}).
$$
We then prove the uniform convergence to $U_{c_*}(r-{\mathrm{ln}}\alpha^\infty(\Theta))$ by examining the difference
$$
\tilde v(t,r,\Theta)=\big\vert v(t,r,\Theta)-e^{r}U_{c_*}(r+\sigma^\infty(t,\Theta))\big\vert
$$
in the region $\{r<t^\delta\}$. For $N=1$, it turned out in \cite{NRR-Brezis} that $\tilde v(t,x)$ was a subsolution of (a perturbation of) the
heat equation
\begin{equation}
\label{e2.301}
\begin{array}{rll}
V_t= &V_{xx} + O(t^{1-\delta}) \, ,& \quad t>0 \, ,\, -t^\delta<x<t^\delta\\
V(t,-t^\delta)=&e^{-t^\delta} \, ,& \quad t>0\\
V(t,t^\delta)=&0 \, ,& \quad t>0.
\end{array}
\end{equation}
The condition at $x=-t^\delta$ simply comes from the fact that $v(t,x)$ decays, by definition, like $e^x$ at $-\infty$. Although the domain might look very large, its first Dirichlet eigenvalue is of the order $t^{-2\delta}$, hence a much larger quantity than the right hand side of \eqref{e2.301}. Thus $V(t,x)$ could be proved to go to 0 uniformly in $x$ as $t\to+\infty$, which implied the sought for convergence result. The same idea will work here again, up to the caveat that $\alpha(\Theta)$ is only Lipschitz in $\Theta$, something that does not go very well with taking a Laplace-Beltrami operator. A regularisation argument, together with some addtional technicalities, will settle the issue.

\noindent Our experience with working with multi-dimensional reaction-diffusion equations is that the main additional difficulty is the transverse diffusion, which, in a very paradoxical way, does not help. This is not a rhetorical argument: its presence is really what prevented convergence in the earlier paper \cite{JMR_VRM2}. This explains why we have to be extra careful with the estimates.

\section{Convergence in the diffusive zone}
\label{s3}

As announced at the end of Section \ref{s2}, the main effort will be to prove a sufficiently strong compactness property for the solutions of \eqref{shift}. Once this is done, convergence to a solution of  the Dirichlet heat equation in self-similar variables, up to a multiplicative coefficient depending on $\Theta$, will follow. In more precise terms, let $\phi_0(\xi)=\xi e^{-\xi^2/4}$, it solves
\begin{equation}
\label{e3.1}
{\cal L}\phi=0\ (\xi>0),\quad\phi(0)=\phi(+\infty)=0.
\end{equation}
Any solution of \eqref{e3.1} is a multiple of $\phi_0$. The main result of this section is the following.
\begin{thm}
\label{t3.1}
Let $w(\tau,\xi)$ be the solution of \eqref{shift} with compactly supported initial datum $w_0$. There exists a Lipschitz function $\Theta\mapsto\alpha(\Theta)$, positive on the unit sphere, such that
$$
\lim_{\tau\to+\infty}w(\tau,\xi,\Theta)=\alpha(\Theta)\phi_0(\xi),
$$
uniformly in $\xi\in\RR_+$ and $\Theta \in \mathbb{S}^{N-1}$.
\end{thm}

\subsection{Sub and super solutions in the diffusive zone, compactness in the $(\tau,\xi)$ variables}

This part is the most technical of the paper, but, many computations being in the spirit of those of 
\cite{NRR-Brezis} or \cite{JMR_VRM2}, we will not detail all of them, rather give their main steps. Let us set
\begin{equation}
\label{e3.100}
\xi_\delta^\pm(\tau)=\pm e^{-(\frac12-\delta)\tau},
\end{equation}
we will often use the notation $\xi_\delta$ and not mention the dependence in $\tau$, as things will - hopefully - be clear from the context. The constant $\delta>0$ will be suitably small and, in any case, less that $1/4$. This will not be made more explicit in the subsequent computations.
The point $\xi=\xi_\delta^+(\tau)$ corresponds, in the $(t,r,\Theta)$ variables, to $r=t^{\delta}$ in the moving frame, that is, far ahead of the supposed location of the front ($r = O(1)$), but not quite as far as the diffusive zone ($r \sim \sqrt{t}$). The point $\xi=\xi_\delta^-(\tau)$ therefore corresponds to $r=-t^{\delta}$, that is, far at the back of the front location, but, again, not quite as far as $-\sqrt t$. 

The main step in this section will be the construction of a subsolution estimating $w(\tau,\xi,\Theta)$ in the region $\{\xi\geq\xi_\delta^+(\tau)\}$ (that is, $r$ starting far ahead of the front) and a super-solution estimating $w(\tau,\xi,\Theta)$ in the region $\{\xi\geq \xi_\delta^-(\tau)\}$, (that is, $r$ starting far at the back of the front). In the $(\tau,\xi)$ variables, the two end points $\xi_\delta^\pm(\tau)$ will rejoin at $\xi=0$ as $\tau\to+\infty$: this will provide an estimate of the solution in the self-similar variables at $\xi\sim0$, whereas the main body of the sub and super solutions will estimate $w$ in the diffusive zone.

For $a>0$, let $\lambda_1(a)$ be the first eigenvalue of the Dirichlet Laplacian on $(-a,a) \subset \RR$ whose eigenfunction $\phi_{1,a}$ has a maximum equal to 1. Note that the maximum is attained at $\xi=0$ and that $\lambda_1(a)=\di\frac{\pi^2}{4a^2}$. For every $a_0>0$ such that $\lambda_1(a_0)\geq100$, we have $a_0<1$. These two easy estimates will be used repeatedly in the course of the section. The real number $a_0$ will be, from then on,  chosen this way. It may also be suitably decreased, independently of all other coefficients and variables. We will set
$$
\lambda_1=\lambda_1(a_0),\ \phi_1(\xi)=\phi_{1,a_0}(\xi).
$$
Let us state the result that will be of use to us in the next section.
\begin{prop}
\label{p3.1}
\textbf{1. Control of $w$ from above and from the back of the front}. \\
There is a pair of positive functions $(q_+(\tau),\zeta_+(\tau))$ that have the following properties.
\begin{enumerate}
\item $\zeta_+$ is bounded and bounded away from 0 by constants that depend only on the initial datum and the constants appearing in the equation, whereas there is $\mu>0$ such that $q_+(\tau)=O(e^{-\mu\tau})$ as $\tau\to+\infty$.
\item For $\xi\geq \xi_\delta^-(\tau)$, we have 
\begin{equation}
\label{e3.101}
w(\tau,\xi,\Theta)\leq \biggl( q_+(\tau)\biggl(\phi_1(\xi-\xi_\delta^-(\tau))+e^{-(\xi-\xi_\delta^-)^2/16}\biggl)+\zeta_+(\tau)\phi_0(\xi-\xi_\delta^-(\tau)) \biggl) e^{-(\xi-\xi_\delta^-)^2/8}.
\end{equation}
\end{enumerate}
 \textbf{2. Control of $w$ from below and from the head of the front}.\\ 
 There is a pair of positive functions $(q_-(\tau),\zeta_-(\tau))$ that satisfy the same estimates as for $q_+$ and $\zeta_+$ in item 1 above, and such that, for $\xi\geq\xi_\delta^+(\tau)$, we have 
\begin{equation}
\label{e3.110}
w(\tau,\xi,\Theta)\geq \biggl(-q_-(\tau) e^{-(\xi-\xi_\delta^+)^2/16}+\zeta_-(\tau)\phi_0(\xi-\xi_\delta^+(\tau)) \biggl) e^{-(\xi-\xi_\delta^+)^2/8}.
\end{equation}
\end{prop}
To prove this proposition, we will make an intermediate step in proving a simple estimate on a particular integro-differential equation.
 Indeed, proving \eqref{e3.101} and \eqref{e3.110} will imply looking for multiple functions of the type $q_\pm$ or $\zeta_\pm$, that will all boil down to solving 
the problem \eqref{e3.102} below. So, 
consider, for any constants $a>0$, $b>0$ and $C>0$, the integro-differential inequality
\begin{equation}
\label{e3.102}
\begin{array}{rll}
\dot q+(a-Ce^{-b\tau})q\leq Ce^{-b\tau}(1+\di\int_0^\tau q(\sigma)d\sigma) & \tau >0\\
q(0)=q_0. &
\end{array}
\end{equation}
with the initial condition $q_0>0$. We have the
\begin{lm}
\label{l3.1}
Equation \eqref{e3.102} has a unique solution $q(\tau)>0$. Moreover there is $K>0$, depending only on $a$, $b$, $C$ and $q_0$ such that 
$$
q(\tau)\leq Ke^{-\frac12 \min(a,b)\tau}.
$$
\end{lm}
\noindent{\bf Proof.} Inequation \eqref{e3.102} is rewritten as
\begin{equation}
\label{e3.103}
\begin{array}{rll}
\dot q+aq\leq Ce^{-b\tau}(1+q+\di\int_0^\tau q(\sigma)d\sigma)&\\
q(0)=q_0.&
\end{array}
\end{equation}
Now, it is enough to prove the existence of $K>0$, depending only on $q_0$ and the coefficients of the problem, such that 
$q(\tau)\leq K$. Indeed, once this is at hand, inequation \eqref{e3.103} implies the inequality
$$
\dot q+aq\leq C e^{-b\tau}(1+K+K\tau),
$$
and the exponential decay of $q$ is deduced from Gronwall's lemma. So, let us concentrate on the global upper bound. 

We first prove that $q$ grows at most exponentially fast: there exists $\Lambda >0$ large enough depending algebraically on $C$ and $q_0$ such that $q(\tau) \leq 2q_0 e^{\Lambda \tau}$ for any $\tau \geq 0$. Indeed, if this is not true, define $\tau_0$ the supremum of $\{\tau \geq 0 \, | \, \forall s\in [0,\tau] \, , \, q(s)\leq 2q_0 e^{\Lambda s}\}$. Then, $\tau_0>0$, $q(\tau_0)=2q_0e^{\Lambda \tau_0}$ and 
$$\dot{q}(\tau_0) \geq  \frac{d}{d\tau} \left(2q_0 e^{\Lambda \tau} \right)_{\tau =\tau_0} = \Lambda q(\tau_0)$$
which implies that $\Lambda^2-C(1+1/q_0)\Lambda -C \leq 0$. This is a contradiction for $\Lambda$ large enough, depending on $q_0$ and $C$.

Let us now improve this exponential bound to a constant. Define $\tau_m>0$ large enough such that
\begin{equation}
\label{e3.105}
\frac32 \frac{C}{a}e^{-b\tau_m}\leq \frac{q_0}2 
\end{equation}
and
\begin{equation}
\label{e3.106}
\hbox{for all }\tau\geq\tau_m \, , \quad  \frac{Ce^{-b\tau}}a(1+\tau)\leq\frac13.
\end{equation}
Assume the existence of $\hat\tau>\tau_m$ such that $q$ has a global maximum over $[0,\hat\tau]$ at $\tau=\hat\tau$. Call $\hat M \geq q_0$ this maximum; at $\tau=\hat\tau$ we have $\dot q\geq0$, so that \eqref{e3.103} becomes, at $\tau=\hat\tau$:
$$
a\hat M\leq Ce^{-b\hat\tau}(1+\hat M+\hat\tau\hat M),
$$
and, because of \eqref{e3.106}:
$$
\hat M\leq\frac32 \frac{C}{a}e^{-b\hat\tau}.
$$
But, thanks to \eqref{e3.105}, we have $\hat M\leq\di\frac{q_0}2$ and a contradiction. Therefore, $\displaystyle \max_{\tau \geq 0} q(\tau)$ is reached over  the interval $[0,\tau_m]$ and we have
$$
\hbox{for all $\tau\geq0$,}\ q(\tau)\leq q(\tau_m) \leq 2q_0 \, e^{\Lambda\tau_m}.
$$
Because $\tau_m$ depends logarithmically on $q_0$, this upper bound grows algebraically in $q_0$. This concludes the proof of lemma \ref{l3.1}. \hfill$\Box$

Before bounding $w(\tau,\xi,\Theta)$ from above and below, we unfortunately need some additional transformations that we detail now. 
The first one is a simple translation in $\xi$ setting the origin at $\xi=\xi_\delta^\pm(\tau)$, and the subsequent one transforms ${\cal L}$ into a self-adjoint operator. This amounts to setting
\begin{equation}
\label{e3.108}
w(\tau,\xi,\Theta)=e^{-(\xi-\xi_\delta^\pm)^2/8}\hat w^\pm(\tau,\hat\xi^\pm,\Theta),\ \hat\xi^\pm=\xi-\xi_\delta^\pm(\tau).
\end{equation}
In order to keep the notations as light as possible, we simply rename the new variable $\hat \xi^\pm$ and the new unknown $\hat w^\pm$ as
$$
\xi:=\xi-\xi_\delta^\pm(\tau),\ \hat w^\pm(\tau,\hat\xi^\pm,\Theta):=w(\tau,\xi,\Theta).
$$
The new equation for $w$ is
\begin{equation}
\label{e3.111}
\left\{
\begin{array}{rll}
\partial_\tau w+{\cal M}w=&l_1(\tau,\xi)\partial_\xi w+l_2(\tau,\xi)w+\di\frac{\Delta_{\Theta}w}{\left(\xi+\xi_\delta^\pm+2e^{\frac{\tau}{2}}-k\tau e^{-\frac{\tau}{2}}\right)^2}-e^{\frac{3\tau}{2}-\frac{\xi^2}{8}-(\xi+\xi_\delta^\pm)e^{\frac{\tau}{2}}}w^2\\
w(\tau,0,\Theta)=&O(e^{\-\mp\delta e^{\tau/2}})\\
w(0,\xi,\Theta)=&w_0(\xi,\Theta)\ \hbox{compactly supported.}
\end{array}
\right.
\end{equation}
The functions $l_1$ and $l_2$ depend on $h$ and are given by
$$l_1(\tau, \xi)=(\frac12-\delta)\xi_\delta^\pm +h(\tau, \xi+\xi_\delta^\pm)e^{-\tau/2} \mbox{ and } l_2(\tau,\xi)=-\frac32-h(\tau, \xi+\xi_\delta^\pm)-\frac{\xi}{4} l_1(\tau,\xi)
$$ 
They satisfy, for  $\xi\geq0$, the following estimates
\begin{equation}
\label{e3.112}
\begin{array}{rll}
\vert l_1(\tau,\xi)\vert\leq &C\biggl(1+\xi\un_{\xi\leq e^{(\frac12-\delta)\tau}}\biggl)e^{-(\frac12-\delta)\tau}\\
\vert l_2(\tau,\xi)\vert\leq &C\biggl(\xi e^{-(\frac12-\delta)\tau}+\un_{\xi\geq e^{(\frac12-\delta)\tau}}\biggl),
\end{array}
\end{equation}
the constant $C$ being universal. The operator ${\cal M}$ writes
\begin{equation}
\label{e3.114}
{\cal M}w=-\partial_{\xi\xi}w+(\frac{\xi^2}{16}-\frac34)w.
\end{equation}
Notice that, if we set
\begin{equation}
\label{e3.115}
\varphi_0(\xi)=\xi e^{-\xi^2/8}=e^{\xi^2/8}\phi_0(\xi),
\end{equation}
we have ${\cal M}\varphi_0=0$. 
We also notice that the Dirichlet condition for $w$ in \eqref{e3.111} is doubly exponentially small if we have taken the origin at $\xi=\xi_\delta^-(\tau)$, whereas it seems doubly exponentially large in the reverse case. This is simply due to the change of variables in item 3 of section 2: $u=e^{-r}v$. So, there is nothing extraordinary in this estimate, at least in the case when we take the origin at the left. When we take the origin at the right, we simply get an over-estimated size that will be of no use to us, because we will seek a bound from below for $w$. There we will just use its positivity.

We are now in a position to bound the solution $w(\tau,\xi,\theta)$ of \eqref{shift} by bounding the new (pair of) solution(s) of \eqref{e3.111}. This is where we will construct a pair of sub and super-solutions. We will concentrate on the upper bound, this will be the subject of the next sub-section. In a subsequent one, we will indicate the modifications necessary for the construction of a subsolution.

\subsubsection{Super-solution}
\label{s3.1.1}

Let $\gamma_1(\xi)$ be a nonnegative smooth function, equal to 1 if $\xi\leq \di\frac{a_0}2$, and zero if $\xi\geq\frac{2a_0}3$, and $\gamma_2(\xi)$ be a smooth nonnegative function, equal to 1 if $\xi\geq 2$ and zero if $\xi\leq1$. 
We will bound the solution $w(\tau,\xi,\Theta)$ from above by a super-solution of the form:
\begin{equation}
\label{e3.116}
\overline w(\tau,\xi,\Theta)=q_1(\tau)\phi_1(\xi)\gamma_1(\xi)+q_2(\tau)\gamma_2(\xi)e^{-\mu\xi^2}+\zeta(\tau)\varphi_0(\xi),
\end{equation}
where $\mu$ is any number in $(0,\di\frac18)$.
Let us set
$$
{\cal N}w=\partial_\tau w+{\cal M}w-l_1(\tau,\xi)\partial_\xi w-l_2(\tau,\xi)w-\di\frac{\Delta_{\Theta}w}{\left(\xi+\xi_\delta^-+2e^{\frac{\tau}{2}}-k\tau e^{-\frac{\tau}{2}}\right)^2}+e^{\frac{3\tau}{2}-\frac{\xi^2}{8}-(\xi+\xi_\delta^-)e^{\frac{\tau}{2}}}w^2.
$$
We want $\overline w$, defined by \eqref{e3.116}, to satisfy ${\cal N}\overline w\geq0$ for $\tau>0,$ $\xi>0$ and $\Theta$ on the unit sphere. A sufficient condition for that is
\begin{equation}
\label{e3.117}
\overline{\cal N}\overline w\geq0,
\end{equation}
with
\begin{equation}
\label{e3.140}
\begin{array}{rll}
\overline{\cal N}w=&\partial_\tau w+{\cal M}w-l_1(\tau,\xi)\partial_\xi w-l_2(\tau,\xi)w\\
&-\di\frac{\Delta_{\Theta}w}{\left(\xi+\xi_\delta^-+2e^{\frac{\tau}{2}}-k\tau e^{-\frac{\tau}{2}}\right)^2};
\end{array}
\end{equation}
in other words we have dropped the positive nonlinear term. We are now looking for sufficient conditions on $q_1$, $q_2$ and $\zeta$ to have $\overline{\cal N}\overline w\geq0$. Throughout the computations, we look for $\zeta$ satisfying moreover
\begin{equation}
\label{e3.118}
\dot\zeta(\tau)\geq0.
\end{equation}
We are going to assume that it is satisfied, then check that it is indeed true. We remark here that we are in the spirit of Fife-McLeod \cite{FML}: because the null space of ${\cal M}$ is not empty, the best we can do with a bare hand computation is estimating the solution, but not prove its convergence, as we have no idea of what multiple of $\varphi_0$ will be picked eventually. In \cite{FML}, a similar computation estimated the position of the front, but did not prove convergence to a wave, as the translation invariance would not permit to guess what translate of the wave would be eventually picked up.

\noindent{\bf 1. The region $0\leq\xi\leq\di\frac{a_0}2$.} 	
This is, in comparison to \cite{NRR-Brezis} and \cite{JMR_VRM2}, the newest part. Here, we have $\gamma_2=0$ since $a_0<1$, so that, using ${\cal M}\varphi_0=0$ and $-\varphi_1''=\lambda_1\phi_1$, we have:
$$
\begin{array}{rll}
\overline{\cal N}\overline w=&\biggl(\dot q_1+(\lambda_1+\di\frac{\xi^2}{16}-\di\frac34)q_1\biggl)\phi_1-\biggl(l_1(\tau,\xi)\phi_1'+l_2(\tau,\xi)\phi_1\biggl)q_1\\
&+\dot\zeta\varphi_0-\zeta\biggl(l_1(\tau,\xi)\varphi_0'+l_2(\tau,\xi)\varphi_0\biggl).
\end{array}
$$
Recall that
$$
\phi_1(\xi)=\cos \left(\frac\pi{2a_0}\xi \right),
$$
so that
$$
\vert\phi_1'(\xi)\vert\leq\sqrt {\lambda_1} \, \phi_1(\xi)\ \hbox{on $[0,\di\frac{a_0}2]$}.
$$
On $[0,a_0/2]$, the functions  $\varphi_0$ and $\varphi_0'$ are estimated by a constant $C$ that - and this is the main thing - does not depend on $a_0$, whereas $\phi_1$ stays above $\sqrt 2/2$ on the same interval. The function $l_1$ is estimated by $e^{-(\frac12-\delta)\tau}$. In the range that we consider, the indicator function appearing in the estimate \eqref{e3.112} of $l_2$ vanishes after a (controlled) finite time, thus we may also estimate $l_2$ by $e^{-(\frac12-\delta)\tau}$. This allows us to assert the existence of $C>0$ independent of $a_0$ such that, if assumption \eqref{e3.118} is true, we have
$$
\frac{\overline{\cal N}\overline w}{\phi_1}\geq\dot q_1+(\lambda_1-C\sqrt\lambda_1-C)q_1-C\zeta e^{-(\frac12-\delta)\tau}.
$$
We choose $a_0>0$ such that $\lambda_1$ is large enough and
$$
\lambda_1-C\sqrt\lambda_1-C\geq1,
$$
this will fix $a_0$ once and for all. And so, a sufficient condition to have $\overline{\cal N}\overline w\geq0$ in this region is
\begin{equation}
\label{e3.119}
\dot q_1+q_1\geq C\zeta e^{-(\frac12-\delta)\tau}.
\end{equation}
{\bf 2. The region $\xi$ large.} 
By this, we mean that $\xi$ will be larger than a constant  $\xi_0 \geq 2$ that we will fix in the course of this section. In any case we have $\gamma_1(\xi)=0$ and $\gamma_2(\xi)=1$. And so,
$$
\begin{array}{rll}
\overline{\cal N}\overline w=&\biggl(\dot q_2+\bigl((\di\frac{1}{16}-4\mu^2)\xi^2+2\mu-\di\frac34\bigl)q_2\biggl)e^{-\mu\xi^2}-\biggl(2\mu l_1(\tau,\xi)\xi-l_2(\tau,\xi)\biggl)q_2e^{-\mu\xi^2}\\
&+\dot\zeta\varphi_0-\zeta\biggl(l_1(\tau,\xi)\varphi_0'+l_2(\tau,\xi)\varphi_0\biggl).
\end{array}
$$
Choose $\mu=\di\frac1{16}$, by  assumption \eqref{e3.118}, we have
$$
\begin{array}{rll}
e^{\xi^2/16}\overline{\cal N}\overline w\geq&\dot q_2+\biggl(\di\frac{3}{64}\xi^2-\frac{\vert l_1(\tau,\xi)\vert\xi}{8}-\vert l_2(\tau,\xi)\vert-\di\frac{5}{8}\biggl)q_2\\
&-\zeta\biggl(l_1(\tau,\xi)\varphi_0'+l_2(\tau,\xi)\varphi_0\biggl) e^{\xi^2/16}.
\end{array}
$$
We estimate $l_i(\tau,\xi)$ as
$$
\vert l_1(\tau,\xi)\vert\leq C \mbox{ and } \vert l_2(\tau, \xi)\vert \leq C(\xi+1) \, .
$$
Thus, the term in factor of $q_2$ can be bounded from below by $\frac{3}{64}\xi^2-\frac98 C \xi -(C+\frac58) $.  Now, we fix $\xi_0$ large enough so that
$$
\frac{3}{64}\xi_0^2-\frac98 C \, \xi_0 -(C+\frac58) \geq 1  \, .
$$
Finally, the function $\varphi_0'$ decays as $\xi^2e^{-\xi^2/8}$ and $\varphi_0(\xi)=\xi e^{-\xi^2/8}$, so that 
$$
\left|l_1(\tau,\xi)\varphi_0'(\xi)+l_2(\tau,\xi)\varphi_0(\xi) \right|  e^{\xi^2/16} \leq C e^{-(\frac12-\delta)\tau}.
$$
Then, the condition $\overline{\cal N}\overline w\geq0$ is satisfied if we have the sufficient condition
\begin{equation}
\label{e3.120}
\dot q_2+q_2\geq C\zeta e^{-(\frac12-\delta)\tau}.
\end{equation}
{\bf 3. The region $\di\frac{a_0}2\leq\xi\leq\xi_0$.} 
Notice that, in this range, the functions $l_i$ may be estimated by ${e^{-(\frac12-\delta)\tau}}$. The functions $\gamma_i$ may take all values between 0 and 1, and their derivatives are bounded. Thus we have
$$
\begin{array}{rll}
\overline{\cal N}\overline w=&\dot\zeta\varphi_0-\zeta\biggl(l_1(\tau,\xi)\varphi_0'+l_2(\tau,\xi)\varphi_0\biggl)\\
&+\dot q_1\gamma_1\phi_1+\dot q_2\gamma_2 e^{-\xi^2/16}+\biggl({\cal M}-l_1(\tau,\xi)\partial_\xi-l_2(\tau,\xi)\biggl)(q_1\gamma_1\phi_1+q_2\gamma_2e^{-\xi^2/16})\\
\geq&\dot\zeta\varphi_0-\zeta\biggl(l_1(\tau,\xi)\varphi_0'+l_2(\tau,\xi)\varphi_0\biggl)+\dot q_1\gamma_1\phi_1+\dot q_2\gamma_2 e^{-\xi^2/16} -C(q_1+q_2).
\end{array}
$$
To render $\overline{\cal N}\overline w$ nonnegative in this range, a sufficient condition is to assume that \eqref{e3.119} and \eqref{e3.120} are satisfied, so that
$\dot q_1\gamma_1\phi_1\geq -q_1$, $\dot q_2\gamma_2e^{-\xi^2/16}\geq -q_2$. Moreover, $\varphi_0$ is bounded away from 0 in our range, so that the final sufficient condition is:
\begin{equation}
\label{e3.121}
\dot\zeta\geq C(q_1+q_2+\zeta e^{-(\frac12-\delta)\tau}).
\end{equation}
It now remains to prove that such functions $q_i$ and $\zeta$ do exist. Setting $q=q_1+q_2$, a sufficient condition for inequalities \eqref{e3.119}, \eqref{e3.120} and \eqref{e3.121} to hold is that the couple $(q,\zeta)$ solves the system
\begin{equation}
\label{e3.122}
\left\{
\begin{array}{rll}
\dot q+q=&C\zeta e^{-(\frac12-\delta)\tau}\\
\dot\zeta=&C(q+\zeta e^{-(\frac12-\delta)\tau}),
\end{array}
\right.
\end{equation}
with suitably large initial data $(q_0,\zeta_0)$. System \eqref{e3.122} has a unique solution by elementary Cauchy theory, and a standard argument shows that $q$ and $\zeta$ are positive throughout their evolution. And this also entails $\dot\zeta\geq0$.

\noindent{\bf Proof of Proposition \ref{p3.1}, Point 1.} 

\noindent
Let $w$ be a solution of \eqref{shift}. Perform transformations \eqref{e3.108} with $\xi^\pm_\delta(\tau)=\xi_\delta^-(\tau)$ so that $w$ is now a solution to \eqref{e3.111}. Define $q$ and $\zeta$ as in \eqref{e3.122} with large initial data $(q_0, \zeta_0)$. Lemma \ref{l3.1} yields the exponential decay for $q$. So, if we choose $q_1=q_2=\di\frac{q}2$ and $q_+=q_1$, 
$\zeta_+=\zeta$, the function $\overline w(\tau,\xi)$, defined in \eqref{e3.116}, is a super-solution to \eqref{e3.111}.
If we choose $q_0$ and $\zeta_0$ large enough, we have $\overline w(0,\xi)\geq w_0(\xi,\Theta)$ for any $\Theta$ on the sphere.  Because $\zeta\geq0$, we have $q(\tau)\geq q_0e^{-\tau}$. Thus
we have, at the expense of choosing $q_0$ even larger:
$$
\overline w(\tau,0)\geq w(\tau,0,\Theta),
$$
the last quantity being less than a double exponential. This finishes the proof. \hfill$\Box$

\subsubsection{Subsolution}

Let us consider again $\gamma_2(\xi)$, a smooth nonnegative function, equal to 1 if $\xi\geq 2$ and zero if $\xi\leq1$. 
This time, we want (instructed by the previous section) to control $w(\tau,\xi,\Theta)$ (after transformations \eqref{e3.108} with $\xi^\pm_\delta(\tau)=\xi_\delta^+(\tau)$) from below by a subsolution of the form:
\begin{equation}
\label{e3.123}
\underline w(\tau,\xi,\Theta)=-q_2(\tau)\gamma_2(\xi)e^{-\xi^2/16}+\zeta(\tau)\varphi_0(\xi),
\end{equation} 
with $q_2>0$. We have dropped the term $q_1\gamma_1\phi_1$ because the function $\underline w$, as just defined, is negative at $\xi=0$, provided $q_2(\tau)\geq0$; this will certainly control $w$ from below at $\xi=0$. 
This time, the nonlinear operator ${\cal N}$ is
$$
\begin{array}{rll}
{\cal N}w=&
\partial_\tau w+{\cal M}w-l_1(\tau,\xi)\partial_\xi w-l_2(\tau,\xi)w\\
&-\di\frac{\Delta_{\Theta}w}{\left(\xi+\xi_\delta^++2e^{\frac{\tau}{2}}-k\tau e^{-\frac{\tau}{2}}\right)^2}+e^{\frac{3\tau}{2}-\frac{\xi^2}{8}-(\xi+\xi_\delta^+)e^{\frac{\tau}{2}}}w^2.
\end{array}
$$
We want $\underline w$, defined by \eqref{e3.123}, to satisfy ${\cal N}\underline w\leq0$ for $\tau>0,$ $\xi>0$ and $\Theta$ on the unit sphere. The $w^2$ term may look bothering; however Point 1 of Proposition \ref{p3.1} is now proved, so that we have, using (in quite a non-optimal way) the proposition:
$$
w^2\leq Cw,
$$
$C$ once again possibly huge. Let us also notice that, for $\xi\geq0$, we have 
$$
\frac{3\tau}2-(\xi+\xi_\delta^+)e^{\tau/2}\leq\frac{3\tau}2-\xi_\delta^+e^{\tau/2}=\frac{3\tau}2-e^{\delta\tau},
$$
so that, all in all, we have for $\tau >0$ and $\xi \geq 0$,
$$
e^{\frac{3\tau}{2}-\frac{\xi^2}{8}-(\xi+\xi_\delta^+)e^{\frac{\tau}{2}}}w^2 \leq Ce^{-(\frac12-\delta)\tau}w.
$$
This term may therefore be included in $l_2(\tau,\xi)$, and 
a sufficient condition for ${\cal N}\underline w\leq0$ is
\begin{equation}
\label{e3.127}
\underline{\cal N} \, \underline w\leq0,
\end{equation}
with
$$
\underline{\cal N}w=\partial_\tau w+{\cal M}w-l_1(\tau,\xi)\partial_\xi w-l_2(\tau,\xi)w-\di\frac{\Delta_{\Theta}w}{\left(\xi+\xi_\delta^++2e^{\frac{\tau}{2}}-k\tau e^{-\frac{\tau}{2}}\right)^2};
$$
the term $l_2$ now incorporating an additional $e^{-(\frac12-\delta)\tau}$. From then on, the computation proceeds in a similar fashion as before, which ends the proof of Proposition \ref{p3.1}.

\subsubsection{Conclusion}

Parabolic regularity yields the boundedness of $\partial_\tau w(\tau,\xi,\Theta)$, $\partial_\xi w(\tau,\xi,\Theta)$ and $\partial_{\xi\xi}w(\tau,\xi,\Theta)$ in terms of the supremum of $w$ on the product of $(\tau-1,\tau+1)\times(\xi-1,\xi+1)$ by the unit sphere. Of course the diffusion in $\Theta$ is degenerate, but it suffices to rescale $\Theta$ by the square root of the diffusion at the point under consideration, and drop the useless estimate in $\Theta$. So, we end up with the following corollary:
\begin{cor}
\label{c3.1}
For $\tau\geq1$, $\xi>0$ and $\Theta$ on the unit sphere, we have
\begin{equation}
\label{e3.125}
\vert \partial_\tau w(\tau,\xi,\Theta)\vert+\vert\partial_\xi w(\tau,\xi,\Theta)\vert\leq Ce^{-\xi^2/16}.
\end{equation}
Moreover, there are two constants $0<\underline q\leq\overline q$ such that, for $\tau\geq1$, $\xi\leq1$ and $\Theta$ on the unit sphere we have
\begin{equation}
\label{e3.126}
\underline q(\xi-e^{-(\frac12-\delta)\tau})\leq w(\tau,\xi,\Theta)\leq \overline q(\xi+e^{-(\frac12-\delta)\tau}).
\end{equation}
\end{cor}

\subsection{Convergence to an angle-depending self-similar solution}

This part completes  the preceding one by proving a compactness property for the variable $\Theta$, entailing the compactness of the trajectories $(w(T+\tau,\xi,\Theta))_{T>0}$ in a weighted $L^\infty$ norm. As the asymptotic problem will simply be the heat equation in the variables $(\tau,\xi)$, convergence will follow. So, let us proceed with compactness.
\begin{prop}
\label{p3.3}
If $w(\tau,\xi,\Theta)$ is the solution of \eqref{e3.111} (with $\hat{\xi}^\pm=\hat{\xi}^-$), there is $C>0$, depending only on the data, such that, for $\tau>0$, $\xi\geq0$ and $\Theta$ on the unit sphere, we have
$$
\vert\nabla_\Theta w(\tau,\xi,\Theta)\vert\leq Ce^{-\xi^2/16}.
$$
\end{prop}
\noindent{\bf Proof.} Let $\Theta_i$ be any coordinate on the unit sphere, and
$$
w_i(\tau,\xi,\Theta)=\partial_{\Theta_i}w(\tau,\xi,\Theta).
$$
As there is no dependence with respect to $\Theta$ in the coefficients of \eqref{e3.111}, the equation for $w_i$ is very similar to that for $w$:
\begin{equation}
\label{e3.300}
\left\{
\begin{array}{rll}
\partial_\tau w_i+{\cal M}w_i=&l_1(\tau,\xi)\partial_\xi w_i+l_2(\tau,\xi)w_i\\
&+\di\frac{\Delta_{\Theta}w_i}{\left(\xi+\xi_\delta^-+2e^{\frac{\tau}{2}}-k\tau e^{-\frac{\tau}{2}}\right)^2}-2e^{\frac{3\tau}2-(\xi+\xi_\delta^-)e^{\frac{\tau}2}-\frac{\xi^2}{8}}ww_i\\
w_i(0,\xi,\Theta)=&\partial_{\Theta_i} w_0(\xi,\Theta)\ \hbox{compactly supported.}
\end{array}
\right.
\end{equation}
We then resort to a classical trick: Multiplying the equation for $w_i$ by the sign of $w_i$ and using Kato's inequality, then finally the fact that $w\geq0$, we find out that $\vert w_i\vert$ solves the inequation
$$
\partial_\tau \vert w_i\vert+{\cal M}\vert w_i\vert- l_1(\tau,\xi)\partial_\xi \vert w_i\vert-l_2(\tau,\xi)\vert w_i\vert-\di\frac{\Delta_{\Theta}w_i}{\left(\xi+\xi_\delta^-+2e^{\frac{\tau}{2}}-k\tau e^{-\frac{\tau}{2}}\right)^2}\leq0.
$$
If now $\overline w_{i}(\xi)$ is the supremum of $\vert w_i(0,\xi,.)\vert$ over the unit sphere, then we have $\vert w_i(\tau,\xi,\Theta)\vert\leq \overline w_i(\tau,\xi)$ with
$$
\left\{
\begin{array}{rll}
\partial_\tau \overline w_i+{\cal M}\overline w_i=&-l_1(\tau,\xi)\partial_\xi\overline  w_i-l_2(\tau,\xi)\overline w_i\\
\overline w_i(0,\xi)=&\overline w_i(\xi)\ \hbox{compactly supported.}
\end{array}
\right.
$$
Moreover, parabolic regularity yields, for the solution $u(t,r,\Theta)$ of \eqref{ref mvt}:
$$
\vert\nabla_\Theta u(t,-t^\delta,\Theta)\vert\leq C(1+t);
$$
this translates into
$$
\vert\nabla_\Theta v(t,-t^\delta,\Theta)\vert\leq C(1+t)e^{-t^\delta},
$$
thus 
$
\overline w_i(\tau,\xi_\delta^-)\leq Ce^{\tau-e^{\delta\tau}}.$
Hence, $\overline w_i$ may be controlled by a super-solution similar to that constructed in Section \ref{s3.1.1}, which proves the proposition. \hfill$\Box$


\noindent {\bf Proof of Theorem \ref{t3.1}.} Corollary \ref{c3.1} and Proposition \ref{p3.3} yield the compactness of the trajectory $(w(T+.,.,.))_{T>0}$ in the $L^\infty_{\tau,\xi,\Theta}$ norm, weighted by $e^{\xi^2/16}$. Therefore, there is a function $w^\infty(\tau,\xi,\Theta)$ and a sequence $(T_n)_n$ going to infinity such that 
\begin{equation}
\label{e3.134}
\lim_{n\to+\infty}e^{\xi^2/16}\vert w(T_n+\tau,\xi,\Theta)-w^\infty(\tau,\xi,\Theta)\vert=0,
\end{equation}
the limit being locally uniform in $\tau$, and uniform in $(\xi,\Theta)$. Moreover, $w^\infty$ is Lipschitz in all its variables, and we have $w^\infty(\tau,0,\Theta)=0$.

On the other hand, for any smooth function $\varphi(\Theta)$ over the unit sphere, consider the integral
$$
w_\varphi(\tau,\xi)=\int_{\mathbb{S}^{N-1}} w(\tau,\xi,\Theta)\varphi(\Theta)d\Theta.
$$
The equation for $w_\varphi$ is once again quite similar as the preceding ones:
$$
\left\{
\begin{array}{rll}
\partial_\tau  w_\varphi+{\cal M} w_\varphi=&-l_1(\tau,\xi)\partial_\xi  w_\varphi-l_2(\tau,\xi) w_\varphi-e^{\frac{3\tau}2-(\xi+\xi_\delta^+)e^{\frac{\tau}2}- \frac{\xi^2}{8}}\di\int_{\mathbb{S}^{N-1}} w^2\varphi d\Theta\\
 w_\varphi(0,\xi)=&\di\int_{\mathbb{S}^{N-1}} w_0(\xi,\Theta)\varphi(\Theta)d\Theta\ \hbox{compactly supported.}
\end{array}
\right.
$$
The same type of super-solution as in Section \ref{s3.1.1} may be constructed for $w_\varphi$. Moreover, the same type of subsolution can also be constructed, as 
we may simply estimate $w_\varphi$ by a constant. This yields the compactness of $w_\varphi$ in the weighted $L^\infty$ norm, but $w_\varphi$ additionally satisfies a standard parabolic equation in the $(\tau,\xi)$ variables. Therefore, parabolic estimates hold, and a subsequence of $(w_\varphi(T+.,.))_{T>0}$ converges, locally in $\tau$, and in the weighted $L^\infty$ norm in $\xi$, to a solution $w^\infty_\varphi$ of 
\begin{equation}
\label{e3.132}
\left\{
\begin{array}{rll}
\partial_\tau w_\varphi^\infty+{\cal M}w_\varphi^\infty&=0,\quad \tau\in\RR \, , \quad \xi\geq0\\
w^\infty_\varphi(\tau,0)=&0.
\end{array} 
\right.
\end{equation}
The same argument as in \cite{NRR-Brezis}, Lemma 5.1 yields the convergence of the full trajectory $(w_\varphi(T+.,.))_{T>0}$ to a steady solution of \eqref{e3.132},
namely, a nontrivial multiple of $\varphi_0$, that we name $\alpha_{\varphi}\varphi_0$.

The functional $\varphi\mapsto\alpha_\varphi$ is a nonnegative functional acting on the set of all continuous functions of the unit sphere. On the other hand, 
\eqref{e3.134} yields, for all $\tau\in\RR$:
$$
\alpha_\varphi\varphi_0(\xi)=\int_{\mathbb{S}^{N-1}} w^\infty(\tau,\xi,\Theta)\varphi(\Theta)d\Theta.
$$
This implies the following cascade of facts. First, the function $w^\infty$ does not depend on $\tau$, we call it $w^\infty(\xi,\Theta)$. Second, the functional $\varphi\mapsto\alpha_\varphi$ is linear, so, combined with positivity, it is a measure that we call $\mu$. Third, we have, for all $\xi>0$:
$$
\int_{\mathbb{S}^{N-1}} \varphi(\Theta)d\mu(\Theta)\varphi_0(\xi)=\int_{\mathbb{S}^{N-1}} w^\infty(\xi,\Theta)\varphi(\Theta)d\Theta.
$$
This entails that $\di\frac{w^\infty(\xi,\Theta)}{\phi_0(\xi)}$ does not depend on $\xi$, call it $\alpha^\infty(\Theta)$. So, we have
$$
\int_{\mathbb{S}^{N-1}} \varphi(\Theta)d\mu(\Theta)\varphi_0(\xi)=\int_{\mathbb{S}^{N-1}} \alpha^\infty(\Theta)\varphi(\Theta)d\Theta.
$$ 
so that $\mu$ is absolutely continuous with respect to the Lebesque measure, $d\mu\Theta)=\alpha^\infty(\Theta)d\Theta$. Because $w^\infty$ is Lipschitz in $\Theta$, this, in the end, implies that $\alpha^\infty$ is Lipschitz. \hfill$\Box$

\section{Convergence to the shifted wave}
\label{s4}

Let $w$ be a solution to \eqref{shift} with compactly supported initial datum $w_0$.  
Setting $\xi_\delta(\tau):=\xi^+_\delta(\tau)$, we start with the following proposition.
\begin{prop}
\label{p4.1}
For every $\e>0$, there is $\tau_\e>0$ (possibly depending also on $\delta$) and $\eta_\e>0$ such that, for all $\tau\geq\tau_\e$ and $\xi\in[-\xi^\delta(\tau),\eta_0]$ we have:
$$
(\alpha^\infty(\Theta)-\e)\, (\xi+\xi^\delta(\tau)) \leq w(\tau,\xi,\Theta)\leq(\alpha^\infty(\Theta)+\e) \, (\xi+\xi^\delta(\tau)).
$$
\end{prop}
This proposition is a consequence of mean value theorem, and the following corollary of Theorem \ref{t3.1}.
\begin{cor}
\label{c3.2}
We have
$$
\lim_{\tau\to+\infty}\partial_\xi w(\tau,0,\Theta)=\alpha^\infty(\Theta).
$$
\end{cor}
\noindent {\bf Proof.} It is enough to prove that 
\begin{equation}
\label{e4.500}
\lim_{T\to+\infty}\partial_\xi w(T+\tau,\xi,\Theta)=\alpha^\infty(\Theta)\phi_0'(\xi),
\end{equation}
uniformly in $\tau$ in every compact with centre 0, $\xi\in[-\xi^\delta(T+\tau),1]$, and $\Theta$ on the unit sphere. For that, it is enough to show the equicontinuity of the family $(\partial_\xi w(T+.,.,.))_{T>0}$; this, combined to the convergence of the family $(w(T+.,.,.))_{T>0}$, implies the convergence of the derivatives. Recall the equation for $w_T:=w(T+.,.,.)$:
\begin{equation}
\label{e5.501}
\left\{
\begin{array}{rll}
\partial_\tau w_T+{\cal M}w_T=&l_1(T+\tau,\xi)\partial_\xi w_T+l_2(T+\tau,\xi)w_T\\
&+\di\frac{\Delta_{\Theta}w_T}{\left(\xi+\xi_\delta^\pm+2e^{\frac{T+\tau}{2}}-k(T+\tau e^{-\frac{T+\tau}{2}}\right)^2}-e^{\frac{3(T+\tau)}{2}-\frac{\xi^2}{8}-(\xi+\xi_\delta^\pm)e^{\frac{(T+\tau}{2}}}w^2_T\\
w_T(\tau,-\xi^\delta,\Theta),\ &\partial_\tau w_T(\tau,-\xi^\delta,\Theta),\ \vert\nabla_\Theta w_T(\tau,-\xi^\delta,\Theta)\vert=O(e^{-e^{\delta(T+\tau)/2}}).
\end{array}
\right.
\end{equation}
The values of $w$, $\partial_\tau w$ and $\nabla_\Theta w$ at the boundary $\{\xi=-\xi^\delta(\tau)\}$ are evaluated from the equation $u=e^{-r}v$, which entails the double exponential for $w$. The derivatives in $\tau$ and $\Theta$ imply the multiplication by a factor $e^\tau$ at most, hence the (innocent) sacrifice of $e^{\tau/2}$ in the exponential controlling $w$.  Standard parabolic regularity results (this would involve a rescaling in $\Theta$ by $e^T$ so as to transform \eqref{e5.501} into a uniformly parabolic equation with smooth coefficients) up to the boundary yield the uniform boundedness of $\partial_{\tau\xi}w$ and $\partial_{\xi\xi}w$. It remains to control $\partial_{\Theta_i\xi}w$, for all $i$. For this, it suffices to differentiate \eqref{e5.501} in $\Theta$, the result is displayed in \eqref{e3.300}, and we already know the boundedness of $\partial_{\Theta_i}w_T$. At the boundary $\{\xi=-\xi^\delta\}$, all derivatives of $\partial_{\Theta_i}w_T$ are controlled by a double exponential, so that parabolic regularity is applicable again. This entails the local boundedness of $\partial_\xi\partial_{\Theta_i}w_T$, hence, as announced, the equicontinuity of the whole family $(\partial_\xi w_T)_{T>0}$. \hfill$\Box$ 

Because $\alpha^\infty$ is only Lipschitz, we will need to use regularisations. If $(\rho_\e)_{\e>0}$ is an approximation of the identity on the unit sphere, we set
\begin{equation}
\label{e4.1}
\alpha^\infty_\e(\Theta)=(\rho_\e*\alpha^\infty)(\Theta).
\end{equation}
Because $\alpha^\infty$ is Lipschitz and positive, we have
$
\alpha^\infty_\e-C\e\leq\alpha^\infty\leq\alpha^\infty_\e+C\e.
$
What is eventually going to be useful to us is the following inequality, for a possibly different constant $C$:
\begin{equation}
\label{e4.3}
(\alpha^\infty_\e(\Theta)-C\e)\, \xi \leq w(\tau,\xi,\Theta)\leq(\alpha^\infty_\e(\Theta)+C\e) \, \xi.
\end{equation}
\noindent{\bf Proof of Theorem \ref{thm1}.} 
We revert to the $(t,r,\Theta)$ variables, and to the function $v(t,r,\Theta)$ defined in Section \ref{s2}. Recall that the equation for $v$ is
\begin{equation}
\label{e4.2} 
\partial_t v= \partial_{rr} v + \biggl(\frac{N-1}{r+2t-k{\mathrm{ln}}t}-\frac{k}t\biggl) \left(\partial_r v - v\right)+\frac{\Delta_\Theta v}{(r+2t-k{\mathrm{ln}}t)^2} -e^{-r} v^2.
\end{equation}
Also recall that the initial unknown $u(t,r,\Theta)$ satisfies
$
u(t,r,\Theta)=e^{-r}v(t,r,\Theta).
$
Fix $\delta_0\in(0,1/100)$. We infer from inequalities \eqref{e4.3} that, for all $\e>0$, applied with $\delta=\di\frac{\delta_0}2$ and $r=t^{\delta_0}$ (so that 
$\xi=e^{-(\frac12-\delta_0)\tau}$, there is $t_\e>0$ such that, for $t\geq t_\e$ we have
\begin{equation}
\label{e4.4}
(\alpha^\infty_\e(\Theta)-C\e)t^{\delta_0}\leq v(t,t^{\delta_0},\Theta)\leq(\alpha^\infty_\e(\Theta)+C\e)t^{\delta_0}.
\end{equation}
In the similar spirit as \cite{NRR-Brezis} and \cite{JMR_VRM2}, we define the upper and lower shifts as
\begin{equation}
\label{e4.6}
U_{c_*}(r+s^{\pm}_\e)\biggl\vert_{r=t^{\delta_0}}=e^{-r}(\alpha^\infty_\e\pm C\e)r\biggl\vert_{r=t^{\delta_0}}.
\end{equation}
Note that $s^{\pm}_\e$ are both well-defined and that the implicit functions theorem yields
$$
s^{\pm}_\e(\Theta)=-\ln(\alpha^\infty_\e\pm C\e)+O(\frac1{t^{\delta_0}}),\ \ \partial_ts^{\pm}_\e(\Theta)=O(\frac1{t^{1+\delta_0}}).
$$
Moreover, the $L^\infty$ norm of $\Delta_\Theta s^{\pm}_\e$ is bounded by a constant that may blow up as $\e\to0$. Let us define $v^\pm_\e$ as the solutions of \eqref{e4.2} for $t\geq t_\e$, $r \in(-t^{\delta_0},t^{\delta_0})$, $\Theta$ on the unit sphere, that have $v(t_\e,r,\Theta)$ as datum at $t=t_\e$, and that satisfy the Dirichlet conditions:
\begin{equation}
\label{e4.7}
v^\pm_\e(t,t^{\delta_0},\Theta)=(\alpha^\infty_\e\pm C\e)t^{\delta_0},\quad v_\e^+(t,-t^{\delta_0},\Theta)=e^{-t^{\delta_0}},\ v_\e^-(t,-t^{\delta_0},\Theta)=0,
\end{equation}
for $t\geq t_\e$ we have
$$
v^-_\e(t,r,\Theta)\leq v(t,r,\Theta)\leq v_\e^+(t,r,\Theta).
$$
The last step of the proof is to prove that the functions $v^\pm_\e(t,r,\Theta)$ converge to $e^{-r}U_{c_*}(r+s^{\pm}_\e)$, uniformly in $r$ and $\Theta$ in their domains. Because $\e$ is arbitrary, this will imply the convergence of $v$. We set (see \cite{NRR-Brezis})
$$
V^\pm_\e(t,r,\Theta)=v^\pm_\e(t,r,\Theta)-e^{-r}U_{c_*}(r+s^{\pm}_\e);
$$
we have
$$
\begin{array}{rll}
\partial_t V^\pm_\e=&\partial_{rr} V^\pm_\e + \biggl(\di\frac{N-1}{r+2t-k{\mathrm{ln}}t}-\frac{k}t\biggl) \left(\partial_r V^\pm_\e - V^\pm_\e\right)\\
&+\di\frac{\Delta_\Theta 
V^\pm_\e}{(r+2t-k{\mathrm{ln}}t)^2} -e^{-r} (U_{c_*}+\underline{v})V^\pm_\e+O(\frac1{t^{1-\delta_0}}).
\end{array}
$$
We use one last time the trick consisting in multiplying the equation by the sign of $V^\pm_\e$, then using the Kato inequality, then the positivity of $U_{c_*}+v$. This yields
$$
\partial_t \vert V^\pm_\e\vert=\partial_{rr} \vert V^\pm_\e\vert + \biggl(\di\frac{N-1}{r+2t-k{\mathrm{ln}}t}-\frac{k}t \biggl) \left(\partial_r \vert V^\pm_\e\vert - \vert V^\pm_\e\vert\right)-\di\frac{\Delta_\Theta 
\vert V^\pm_\e\vert}{(r+2t-k{\mathrm{ln}}t)^2}+O(\frac1{t^{1-\delta_0}}).
$$
Thus, $\vert V^\pm_\e(t,r,\Theta)\vert\leq \overline V^\pm_\e(t,r)$, with
$$
\begin{array}{rll}
&\partial_t \overline V^\pm_\e=\partial_{rr} \overline V^\pm_\e+ (\di\frac{N-1}{r+2t-k{\mathrm{ln}}t}-\frac{k}t) \left(\partial_r \overline V^\pm_\e
- \overline V^\pm_\e\right)+O(\di\frac1{t^{1-\delta_0}})\ \ (t\geq t_\e,\  r\in(-t^{\delta_0},t^{\delta_0}))\\
&\\
&\overline{V}^\pm_\e(t,-t^{\delta_0})=e^{-t^{\delta_0}},\ \ \overline{V}_\e^\pm(t,t^{\delta_0})=0\quad \overline{V}^\pm_e(t_\e,r)=0.
\end{array}
$$
Just as in \cite{NRR-Brezis} we infer that both functions $\overline V^\pm_\e(t,.)$ converge to 0 as $t\to+\infty$. The reason is that the equation has lower order coefficients and right handside of order less than $\di\frac1t$, whereas the first eigenvalue of the Dirichlet Laplacian on $(-t^{\delta_0},
t^{\delta_0})$ is of order $t^{-2\delta_0}$\hfill$\Box$

 
\section{Discussion}
Let us first mention that our result remains valid for more general nonlinearities. For an equation of the form
$$
\partial_t u=\Delta u+f(u),\quad t>0,\ x\in\RR^N,
$$
it suffices to assume that $f$ is concave and positive on $(0,1)$, with $f(0)=f(1)=0$. Thus $f'(0)>0$ and the bottom speed is given by $c_*=2\sqrt{f'(0)}$. Our result becomes the existence of a function $s$ defined on the unit sphere such that
$$
u(t,x)=U_{c_*}\biggl(\vert x\vert+c_*t-\frac{N+2}{c_*}{\mathrm {ln}}+s^\infty(\frac{x}{\vert x\vert})\biggl)+o_{t\to+\infty}(1).
$$ 
In the course of the proof, the nonlinear term is no more $u^2$ but $g(u)=f'(0)u-f(u)$,  which is positive and nondecreasing on $(0,1)$. It is not clear to us whether the result would subsist by merely assuming $f(u)\leq f'(0)u$. What would probably be true is a statement of the form
$$
u(t,x)=U_{c_*}\biggl(\vert x\vert+c_*t-\frac{N+2}{c_*}{\mathrm {ln}}t+s^\infty(t,\frac{x}{\vert x\vert})\biggl)+o_{t\to+\infty}(1),
$$
with $s^\infty(t,\Theta)=O(1)$. Let us also mention that we could have given a slightly different version of Theorem \ref{thm1} by stating that, for every direction $e$, $\vert e\vert=1$, then
$$
\{u(t,x)=\lambda\}\cap\{x=re,\ r>0\}\subset\{r=c_*t-\frac{N+2}{c_*}{\mathrm {ln}}t-s^\infty(e)+U_{c_*}^{-1}(\lambda)+o_{t\to+\infty}(1)\}.
$$
The analysis of the solution on the diffusive zone would have been slightly simpler, in the sense that we would not have had to handle an asymptotically degenerate diffusion in $e$. On the other hand, recovering the convergence at the $O(1)$ spatial scale would have been more delicate. Additionally, this would not have proved the Lipschitz regularity of $s$ in $e$. This last approach is, sometimes, better tailored to the geometric situation, where the front has a preferered direction of propagation. This is the case in the forthcoming paper \cite{CR}, where the Fisher-KPP invasion occurs orthogonally to a line of fast diffusion.

We may adapt the preceding ideas to asymptotically homogeneous models of the form
\begin{equation}
\label{e5.1}
\partial_t u =\Delta u +\mu(x)u-u^2,\quad (t>0,\ x\in\RR^N)
\end{equation}
where the function $\nu(x):=\mu-1$ satisfies
$$
\nu(x)=\frac\lambda{\vert x\vert^\alpha}+O_{\vert x\vert\to+\infty}(\frac1{\vert x\vert^{\alpha+\delta}}),\ \vert\nabla\nu(x)\vert=\frac{\alpha\lambda}{\vert x\vert^{1+\alpha}}+O_{\vert x\vert\to+\infty}(\frac1{\vert x\vert^{\alpha+1+\delta}}).
$$
Theorem \ref{thm1} becomes
\begin{thm}
\label{thm2}
Let $u_0$ satisfy assumption \eqref{assumption}. There is a Lipschitz function $s^\infty(\Theta)$, defined on the unit sphere of $\RR^N$, such that the solution $u$ of \eqref{KPP} emanating from $u_0$ satisfies
$$
u(t,x)=U_{c_*}\biggl(\vert x\vert-c_*t+\frac{N+2}{c_*}{\mathrm{ln}}t+s^\infty(\frac{x}{\vert x\vert})\biggl)+o_{t\to+\infty}(1),
$$
if $\alpha>1$, and
$$
u(t,x)=U_{c_*}\biggl(\vert x\vert-c_*t+\frac{N+2-\lambda}{c_*}{\mathrm{ln}}t+s^\infty(\frac{x}{\vert x\vert})\biggl)+o_{t\to+\infty}(1),
$$
if $\alpha=1$.
\end{thm}
The shift $\di\frac{N+2-\lambda}{c_*}$ had already been identified by Ducrot \cite{D}, up to $O(1)$ terms. His assumptions are more general than ours, in the sense that he neither requires the gradient estimate on $\nu$, nor the quantitative estimate for $\nu(x)-\di\frac\lambda{\vert x\vert^\alpha}$.  However, our result goes one step further. Theorem \ref{thm2} would probably hold without the error estimate on $\nu(x)$, one would simply need to be more careful in the construction of sub and super solutions. On the other hand, we have not tried to push the limits of validity of Theorem \ref{thm2}, and this might well be quite an interesting question.

The proof of Theorem \ref{thm2} goes exactly along the same lines as that of Theorem \ref{thm1} for $\alpha>1$, the term $\nu(x)$ being thrown into the perturbative terms $l_i(\tau,\xi)$. Of course they now depend on $\Theta$, but in a smooth and exponentially small in time fashion, so they do not require any additional arguments. When $\alpha=1$, the same algebraic steps as in Section 2 revel the presence of a nonperturbative term in equation \ref{dec onde}. More precisely, this equation becomes
\begin{equation}
\label{e5.2}
\begin{array}{rll} 
\partial_t v= &\partial_{rr} v + \di(\frac{N-1}{r+2t-k{\mathrm{ln}}t}-\frac{k}t) \partial_r v-\di{(\frac{N-1}{r+2t-k{\mathrm{ln}}t}-\frac{k}t-\frac\lambda{c_*t})}v\\
&+\di{\biggl(\nu(r+c*t-k{\mathrm{ln}}t,\Theta)-\frac\lambda{c_*t}\biggl)v+\frac{\Delta_\Theta v}{(r+2t-k{\mathrm{ln}}t)^2}} -e^{-r} v^2.
\end{array}
\end{equation}
To identify $k$ we simply have to make sure that equation \eqref{e5.2} behaves like the Dirichlet heat equation, perturbed by higher order terms; thus the formula \eqref{e2.1000} becomes
$$
\frac{N-1}{c_*}-\frac\lambda{c_*}-k=-\frac32,
$$
hence the shift. The remaining terms will be, in the self-similar variables, exponentially decreasing terms. The method used to prove a gradient estimate in $\Theta$ for $v$ will then work exactly as in Proposition \ref{p3.3}, thanks to the estimate on $\vert\nabla\nu\vert$.

We finally mention that we leave open the question of higher order expansion of the shift, which is also quite an interesting question.

{\footnotesize
 
}

\end{document}